# A procedure for the construction of non-stationary Riccati-type flows for incompressible 3D Navier-Stokes equations

**Sergey V. Ershkov**

Institute for Time Nature Explorations,

M.V. Lomonosov's Moscow State University,

Leninskie gory, 1-12, Moscow 119991, Russia

e-mail: sergej-ershkov@yandex.ru

In fluid mechanics, a lot of authors have been executing their researches to obtain the analytical solutions of Navier-Stokes equations, even for 3D case of *compressible* gas flow or 3D case of *non-stationary* flow of incompressible fluid. But there is an essential deficiency of non-stationary solutions indeed.

We explore the ansatz of derivation of *non-stationary* solution for the Navier-Stokes equations in the case of incompressible flow, which was suggested earlier.

In general case, such a solution should be obtained from the mixed system of 2 Riccati ordinary differential equations (in regard to the time-parameter $t$). But we find an elegant way to simplify it to the proper analytical presentation of exact solution (such a solution is exponentially decreasing to zero for $t$ going to infinity $\infty$). Also it has to be specified that the solutions that are constructed can be considered as a class of perturbation absorbed exponentially as $t$ going to infinity $\infty$ by the null solution.

# 1. Introduction, the Navier-Stokes system of equations.

In accordance with [1-3], the Navier-Stokes system of equations for incompressible flow of Newtonian fluids should be presented in the Cartesian coordinates as below (*under the proper initial conditions*):

$$\nabla \cdot \vec{u} = 0, \qquad (1.1)$$

$$\frac{\partial \vec{u}}{\partial t} + (\vec{u} \cdot \nabla)\vec{u} = -\frac{\nabla p}{\rho} + \nu \cdot \nabla^2 \vec{u} + \vec{F}, \qquad (1.2)$$

- where ***u*** is the flow velocity, a vector field; $\rho$ is the fluid density, $p$ is the pressure, $\nu$ is the kinematic viscosity, and ***F*** represents external force (*per unit of mass in a volume*) acting on the fluid. Let us also choose the *Ox* axis coincides to the main direction of flow propagation. Besides, we assume here external force ***F*** above to be the force, which has a potential $\phi$ represented by ***F*** = -$\nabla \phi$.

As for the domain in which the flow occurs and the boundary conditions, let us consider only the Cauchy problem in the whole space.

We should note that the equation of momentum (1.2) of Navier-Stokes system of equations (1.1)-(1.2) could be presented as a direct sum of 3 equations below [4]:

$$\begin{cases} \dfrac{1}{2}\nabla\{(\vec{u}_p + \vec{u}_w)^2\} + \dfrac{\nabla p}{\rho} + \nabla \phi = 0, & (1.3) \\[1em] \dfrac{\partial \vec{u}_p}{\partial t} = \vec{u}_p \times \vec{w} + \vec{f}, & (1.4) \\[1em] \dfrac{\partial \vec{u}_w}{\partial t} = \nu \cdot \nabla^2 \vec{u}_w, & (1.5) \end{cases}$$

- where $f = \vec{u}_w \times \vec{w}$ ; besides, $\vec{u}_p$ is *an irrotational* (*curl-free*) field of flow velocity, and $\vec{u}_w$ - is *a solenoidal* (*divergence-free*) field of flow velocity which generates a curl field $w$, according to the Helmholtz fundamental theorem of vector calculus [5]:

$$\nabla \times \vec{u} \equiv \vec{w}, \qquad \vec{u} \equiv \vec{u}_p + \vec{u}_w,$$

$$\nabla \cdot \vec{u}_w \equiv 0, \qquad \nabla \times (\vec{u}_p) \equiv 0.$$

So, if we solve Eq. (1.5) in regard to the components of vector $\vec{u}_w$, we could substitute it into the Eq. (1.4) to find a solution for $\vec{u}_p$; then we could substitute the derived components $\vec{u}_w$, $\vec{u}_p$ the appropriate expressions in Eq. (1.3) for the obtaining of the proper expression for vector function $\nabla p$:

$$\frac{\nabla p}{\rho} = -\nabla \phi - \frac{1}{2}\nabla\{(\vec{u}_p + \vec{u}_w)^2\},$$

- here $\phi$ - is the also time-dependent function, in general case.

2. <u>**General presentation of the curl-free field of flow, $\vec{u}_p = \{U, V, W\}$.**</u>

Eq. (1.4) is known to be [5] the system of 3 linear PDEs (in regard to the time-parameter $t$) for 3 unknown functions: $U$, $V$, $W$, where $\vec{u}_p = \{U, V, W\}$.

Besides, (1.4) is the system of 3 linear differential equations with all the coefficients depending on time $t$. In accordance with [5] p.71, the general solution of such a system should be given as below:

$$\chi_p(t) = \sum_{v=1}^{3} \zeta_{v,p} \cdot \left( \int\left(\frac{\Delta_v}{\Delta}\right) dt + C_v \right), \quad (p=1,2,3), \quad \Delta = \begin{vmatrix} \zeta_{1,1} & \zeta_{1,2} & \zeta_{1,3} \\ \zeta_{2,1} & \zeta_{2,2} & \zeta_{2,3} \\ \zeta_{3,1} & \zeta_{3,2} & \zeta_{3,3} \end{vmatrix} \quad (2.1)$$

- where $\{\chi_p\}$ - are the fundamental system of solutions of Eqs. (1.4): $\{U, V, W\}$, in regard to the time-parameter $t$; $\{\zeta_{v,p}\}$ - are the fundamental system of solutions of *the corresponding homogeneous* variant of (1.4), $\{C_v\}$ - are the set of functions, not depending on time $t$; besides, here we note as below:

$$\Delta_1 = \begin{vmatrix} f_x & f_y & f_z \\ \zeta_{2,1} & \zeta_{2,2} & \zeta_{2,3} \\ \zeta_{3,1} & \zeta_{3,2} & \zeta_{3,3} \end{vmatrix}, \quad \Delta_2 = \begin{vmatrix} \zeta_{1,1} & \zeta_{1,2} & \zeta_{1,3} \\ f_x & f_y & f_z \\ \zeta_{3,1} & \zeta_{3,2} & \zeta_{3,3} \end{vmatrix}, \quad \Delta_3 = \begin{vmatrix} \zeta_{1,1} & \zeta_{1,2} & \xi_{1,3} \\ \zeta_{2,1} & \zeta_{2,2} & \zeta_{2,3} \\ f_x & f_y & f_z \end{vmatrix}.$$

It means that the system of Eqs. (1.4) could be considered as having been solved if we obtain a general solution of *the corresponding homogeneous* system (1.4). Let us search for such a general solution of the corresponding homogeneous system (1.4) as below

$$\frac{\partial U}{\partial t} = V \cdot w_z - W \cdot w_y,$$

$$\frac{\partial V}{\partial t} = W \cdot w_x - U \cdot w_z, \qquad (2.2)$$

$$\frac{\partial W}{\partial t} = U \cdot w_y - V \cdot w_x.$$

The system of Eqs. (2.2) has *the analytical* way to present the general solution [5] (in regard to the time-parameter $t$) for 3 unknown functions $U, V, W$:

$$U = \frac{(\gamma - W) \cdot (\xi - \eta^{-1})}{2}, \quad V = -\frac{(\gamma - W) \cdot i \cdot (\xi + \eta^{-1})}{2},$$

$$W = \gamma \cdot \frac{\left(1 + \frac{\eta}{\xi}\right)}{\left(1 - \frac{\eta}{\xi}\right)},$$

(2.3)

- where $\gamma(x, y, z)$ – is some arbitrary (real) function, given by the initial conditions; for auxillary functions $\xi(t)$, $\eta(t)$ of *complex* value we could obtain the appropriate *Riccati* equations as below:

$$\xi' = \left(\frac{w_y + i \cdot w_x}{2}\right) \cdot \xi^2 - i \cdot w_z \cdot \xi + \left(\frac{w_y - i \cdot w_x}{2}\right), \qquad (2.4)$$

$$\eta' = \frac{(w_y + i \cdot w_x)}{2} \cdot \eta^2 - i \cdot w_z \cdot \eta + \frac{(w_y - i \cdot w_x)}{2} \qquad (2.5)$$

- besides, we should note that:

$$\eta^{-1} = -\overline{\xi}, \qquad (*)$$

- that's why all the components $\{U, V, W\}$ (2.3) are *the real* functions in any case.

Also, according to the *continuity* equation (1.1), the appropriate restriction should be valid for identifying of function $\gamma(x, y, z)$ and the set of functions $\{C_v(x, y, z)\}$ in (2.1):

$$\frac{\partial U}{\partial x} + \frac{\partial V}{\partial y} + \frac{\partial W}{\partial z} = 0 \qquad (2.6)$$

- which is the PDE-equation of the 1-st kind; $w_x$, $w_y$, $w_z$ depend on variables $(x,y,z,t)$.

So, the existence of the solution for Navier-Stokes system of equations (2.6), (1.3)-(1.5) (*which is equivalent to the initial system of Navier-Stokes* (1.1)+(1.2) *in the sense of existence and smoothness of a solution*) is proved to be the question of existence of the proper function $\gamma(x, y, z)$ and the set of functions $\{C_v(x, y, z)\}$ in (2.1) of so kind that the PDE-equation (2.6) should be satisfied under the given initial conditions.

Besides, *non-homogeneous* solution (2.1) from presentation of a solution (2.3) should generate a zero curl: $\boldsymbol{u}_P$ is *an irrotational* (*curl-free*) field of flow velocity; it means an additional restrictions at choosing of functions $\gamma(x, y, z)$ and $\{C_v(x, y, z)\}$ in (2.1):

$$\frac{\partial W}{\partial y} - \frac{\partial V}{\partial z} = 0, \quad \frac{\partial U}{\partial z} - \frac{\partial W}{\partial x} = 0, \quad \frac{\partial V}{\partial x} - \frac{\partial U}{\partial y} = 0 \qquad (2.7)$$

But we should consider here only the *helical* type of solutions for the *solenoidal* part of flow velocity $u_w$ (for which is valid $f = u_w \times w = 0$); it means that the vorticity vector, determined from the equation below:

$$\frac{\partial \vec{w}}{\partial t} = \nu \cdot \nabla^2 \vec{w},$$

- is parallel to the velocity vector (1.5) at every point of the fluid. So, it should be sufficient to consider the *homogeneous* solution (2.3) for such a case.

### 3. The PDE-system for the curl-free $u_p = \{U, V, W\}$ part of solution.

Let us present the auxillary functions $\xi(t)$, $\eta(t)$ (2.4)-(2.5) of *complex* value as below

$$\eta = a + b \cdot i, \quad \xi = c + d \cdot i, \tag{3.1}$$

- then due to (*), we obtain the non-linear dependence of coefficients $c$, $d$ on coefficients $a$, $b$:

$$\eta^{-1} = \frac{a}{a^2 + b^2} - \frac{b}{a^2 + b^2} \cdot i = -\overline{\xi} = -(c - d \cdot i),$$

$$\Rightarrow \begin{cases} c = -\left(\dfrac{a}{a^2 + b^2}\right), \\ \\ d \cdot i = -\left(\dfrac{b}{a^2 + b^2}\right) \cdot i. \end{cases} \tag{3.2}$$

So, using (3.1)-(3.2), we could represent the solution (2.3) in other form as below:

$$U = -\gamma \cdot \left(\frac{2a}{1 + (a^2 + b^2)}\right), \quad V = -\gamma \cdot \left(\frac{2b}{1 + (a^2 + b^2)}\right),$$

$$W = \gamma \cdot \left(\frac{1 - (a^2 + b^2)}{1 + (a^2 + b^2)}\right). \tag{3.3}$$

Analyzing the expressions (3.3) above, it is explicitly obvious that we should explore only the function η(*t*) and appropriate dynamics of coefficients *a*(*x,y,z, t*), *b*(*x,y,z, t*).

So, we obtain from the *Riccati* equation (2.5) for η (*t*):

$$\eta' = \frac{(w_y + i \cdot w_x)}{2} \cdot \eta^2 - i \cdot w_z \cdot \eta + \frac{(w_y - i \cdot w_x)}{2}, \quad \Rightarrow$$

$$(a+b\cdot i)' = \frac{(w_y + i \cdot w_x)}{2} \cdot (a^2 + 2a\cdot b\cdot i - b^2) - i\cdot w_z \cdot (a+b\cdot i) + \frac{(w_y - i \cdot w_x)}{2}$$

$$\Rightarrow (a+b\cdot i)' = \left( \frac{w_y}{2} \cdot a^2 - b^2 \frac{w_y}{2} - w_x \cdot a \cdot b + w_z \cdot b + \frac{w_y}{2} \right) + \quad (3.4)$$

$$+ (w_y \cdot a \cdot b + \frac{w_x}{2} \cdot a^2 - \frac{w_x}{2} \cdot b^2 - w_z \cdot a - \frac{w_x}{2}) \cdot i$$

The last ordinary differential equation (3.4) of *complex* value (in regard to the time-parameter *t*) could be presented as a system of equations, which describes appropriate dynamics for the real and imaginary components of function η(*t*):

$$\begin{cases} a' = \frac{w_y}{2} \cdot a^2 - (w_x \cdot b) \cdot a - \left( \frac{w_y}{2}(b^2 - 1) - w_z \cdot b \right), \\ \\ b' \cdot i = -\left( \frac{w_x}{2} \cdot b^2 - (w_y \cdot a) \cdot b - \left( \frac{w_x}{2} \cdot (a^2 - 1) - w_z \cdot a \right) \right) \cdot i . \end{cases} \quad (3.5)$$

### 4. <u>Exact solutions of PDE-system for *u* $_p$ = {*U, V, W*}.</u>

Eqs. (3.5) is the mutual system of 2 *Riccati* ordinary differential equations, which has no analytical solution in general case [5-6]:

$$\begin{cases} a' = \frac{w_y}{2} \cdot a^2 - (w_x \cdot b) \cdot a - \frac{w_y}{2}(b^2 - 1) + w_z \cdot b , \\ \\ b' = -\frac{w_x}{2} \cdot b^2 + (w_y \cdot a) \cdot b + \frac{w_x}{2} \cdot (a^2 - 1) - w_z \cdot a . \end{cases}$$

Indeed, if we multiply the 1-st of Eqs. (3.5) on $w_y$, the 2-nd Eq. on $w_x$, then summarize them one to each other properly, we should obtain:

$$w_y \cdot a' - \frac{1}{2}((w_x)^2 + (w_y)^2) \cdot a^2 + w_x \cdot w_z \cdot a + \frac{(w_x)^2}{2} =$$

$$= -w_x \cdot b' - \frac{1}{2}((w_x)^2 + (w_y)^2) \cdot b^2 + w_y \cdot w_z \cdot b + \frac{(w_y)^2}{2}.$$

(4.1)

Equation (4.1) above is the classical *Riccati* ODE. It describes the evolution of function $a(x,y,z, t)$ in dependence on the function $b(x,y,z, t)$ together with functions $\{w_x, w_y, w_z\}$ in regard to the time-parameter *t*; such a *Riccati* ODE has no analytical solution in general case [5]:

$$a' = A \cdot a^2 + B \cdot a + D.$$

$$A = \frac{1}{2} \frac{((w_x)^2 + (w_y)^2)}{w_y}, \quad B = -\left(\frac{w_x \cdot w_z}{w_y}\right), \quad (4.2)$$

$$D = -\frac{w_x}{w_y} \cdot b' - \frac{1}{2}\frac{((w_x)^2 + (w_y)^2)}{w_y} \cdot b^2 + w_z \cdot b + \frac{w_y}{2} - \frac{(w_x)^2}{2w_y}.$$

But nevertheless, some of the important partial solutions of (4.2) should be considered properly [5]:

1) If an additional condition is valid for the coefficients *A*, *B*, *D* in (4.2), as below

$$\delta^2 \cdot A + \varepsilon \cdot \delta \cdot B + \varepsilon^2 \cdot D = 0, \ \{\varepsilon, \delta\} = \text{const}, \ |\varepsilon| + |\delta| > 0 \quad (4.3)$$

- then, by change of variables $a = \delta \cdot \varepsilon^{-1} + s(t)$, Eq. (4.2) could be reduced to the *Bernoulli* equation:

$$s' = A \cdot s^2 + (2\delta \cdot \varepsilon^{-1} \cdot A + B) \cdot s,$$

- which could be easily solved, so we should obtain as a result that the function *a* depends only on the functions $\{w_x, w_y, w_z\}$ in regard to the time-parameter *t*. Having obtained the appropriate expression for the function *a*, we should substitute it to the 2-nd of Eqs. (3.5) for resolving it in regard to the function *b*. Then, we should restrict the choosing of the appropriate functions $\{w_x, w_y, w_z\}$ in accordance with equality (4.3).

2) If an additional condition is valid for coefficients *A*, *B*, *D* in (4.2) as below

$$4D = (B^2 / A) - 2(B / A)' , \qquad (4.4)$$

- then, solution of (4.2) should be $a = - (B / 2A)$; so, we also should obtain that the function *a* depends only on the functions $\{w_x, w_y, w_z\}$ in regard to the time-parameter *t*.

Again, we should substitute expression for function *a* to the 2-nd of Eqs. (3.5) for resolving it in regard to function *b*. Then, we should restrict the choosing of the appropriate functions $\{w_x, w_y, w_z\}$ in accordance with equality (4.4).

As we can see from examples 1) - 2) above, we really could obtain the partial solutions of Eq. (4.1), but in each case we should restrict the choosing of the appropriate functions $\{w_x, w_y, w_z\}$ in accordance with equalities (4.3) or (4.4).

Let also consider another case of exact solution for system (3.5) than (4.3)-(4.4); namely, let us multiply both parts of the 1-st equation of (3.5) to $a(x,y,z, t)$, but the 2-nd equation of (3.5) to the $b(x,y,z, t)$ properly, then we should sum it one to each other as below

$$\frac{1}{2}(a^2 + b^2)' = \left(\frac{w_y}{2}\right) \cdot a \cdot (a^2 + b^2) - \left(\frac{w_x}{2}\right) \cdot b \cdot (a^2 + b^2) + \left(\frac{w_y}{2}\right) \cdot a - \left(\frac{w_x}{2}\right) \cdot b ,$$

- or

$$\frac{1}{2}(a^2 + b^2 + 1)' = \frac{1}{2}(a^2 + b^2 + 1) \cdot (w_y \cdot a - w_x \cdot b), \qquad (4.5)$$

- where the last Eq. (4.5) could be easily integrated if

$$\left(w_y \cdot a - w_x \cdot b\right) = 0, \quad \Rightarrow \quad b = \sqrt{C^2 - a^2} \qquad (4.6)$$

- here $C$ = const ($|a| \leq C$, $a \neq 0$). It yields from the 1-st of Eqs. (3.5) as below:

$$a' = \frac{w_x}{2} \cdot \frac{b}{a} \cdot a^2 - (w_x \cdot b) \cdot a - \frac{w_x}{2} \cdot \frac{b}{a} \cdot b^2 + w_z \cdot b + \frac{w_x}{2} \cdot \frac{b}{a},$$

- or we should obtain from (4.6)

$$a' = \frac{w_x}{2} \cdot \frac{\sqrt{C^2 - a^2}}{a} \cdot a^2 - (w_x \cdot \sqrt{C^2 - a^2}) \cdot a -$$

$$- \frac{w_x}{2} \cdot \frac{\sqrt{C^2 - a^2}}{a} \cdot (C^2 - a^2) + w_z \cdot \sqrt{C^2 - a^2} + \frac{w_x}{2} \cdot \frac{\sqrt{C^2 - a^2}}{a},$$

$$\Rightarrow \quad a' = (1 - C^2) \cdot \frac{w_x}{2} \cdot \frac{\sqrt{C^2 - a^2}}{a} + w_z \cdot \sqrt{C^2 - a^2} \qquad (4.7)$$

Eq. (4.7) gives us not less than 2 elegant opportunities (for integrating of such an equation):

1) $C = 1$, in such a case function $w_z$ does not depend on functions $w_x$ or $w_y$:

$$\int \frac{da}{\sqrt{1 - a^2}} = \int w_z \, dt$$

$$\Rightarrow \quad a\,(x, y, z, t) = \sin\left(\int w_z \, dt\right) \qquad (4.8)$$

We should note that solution (3.3) with components (4.6), (4.8) does not satisfy to the continuity equation (2.6) along with the demands of zero curl (2.7).

2) $(w_z / w_x) = R\,(a, \sqrt{(C^2 + a^2)})$, where $R$ – is a *rational* function, consisting of polynomial in regard to the function $a$ and, besides, of polynomial in regard to the term $\sqrt{(C^2 + a^2)}$ of extent not more than 2, see for example [7]:

$$R(a, \sqrt{(C^2+a^2)}) \sim a \cdot (\sqrt{(C^2+a^2)})^2 \text{ or } R(a, \sqrt{(C^2+a^2)}) \sim a^3 \cdot \sqrt{(C^2+a^2)}.$$

In this case, the left part of (4.7) could be transformed to the proper *elliptic* integral in regard to the function *a*:

$$\int \frac{da}{\frac{(1-C^2)}{2} \cdot \frac{\sqrt{C^2-a^2}}{a} + R(a,\sqrt{C^2-a^2}) \cdot \sqrt{C^2-a^2}} = \int w_x \, dt \qquad (4.9)$$

It is well-known fact that *elliptical* integral is a generalization of the class of inverse periodic functions [6-7]. Thus, by the proper obtaining of re-inverse dependence of a solution (4.9) from the time-parameter *t* we could present the expression for the function *a* as a set of periodic cycles.

Let us obtain the simple exact solutions of system (3.5). If we assume an appropriate simplifications to be valid as below

$$b = \alpha \cdot a \ (\alpha = const \neq 0), \quad w_x = -\alpha \cdot w_y, \quad w_z = 0, \qquad (4.10)$$

- each of equations of system (3.5) could be reduced properly to the more simple form as below, according to (4.10):

$$a' = \left(\frac{w_y}{2} - w_x \cdot \alpha - \frac{w_y}{2}\alpha^2\right) \cdot a^2 + \frac{w_y}{2} \qquad (4.11)$$

The left part of Eq. (4.11) could be transformed to the appropriate *elliptical* integral [7] {in regard to the function *a*}:

$$\frac{da}{\left(a^2 + \frac{1}{(\alpha^2+1)}\right)} = \left(\frac{(\alpha^2+1) \cdot w_y}{2}\right) dt$$

$$\Rightarrow \arctan((\alpha^2+1) \cdot a) = \frac{(\alpha^2+1)}{2} \cdot \int w_y \, dt, \qquad (4.12)$$

$$a(x, y, z, t) = \frac{\tan\left(\frac{(\alpha^2+1)}{2} \cdot \int w_y \, dt\right)}{(\alpha^2+1)}.$$

At Fig.1 we schematically imagine the plot of time-dependent component $U(t)$ for $t > 0$, in accordance with the expression for appropriate component of solution (3.3), (4.12) and under the *self-similarity* assumption below ($\alpha_y = $ const), see (7.2) for example:

$$w_y = \exp(-\alpha_y^2 \cdot t) \cdot w_y(x, y, z)$$

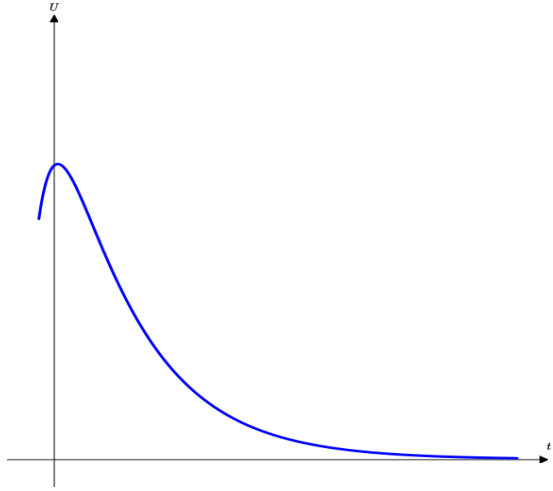

## 5. **Variable curl, helical type of flow.**

We could derive the components of *solenoidal* parts of the solution $\boldsymbol{u}_w = \beta \cdot \boldsymbol{w}$ (1.5), corresponding to the *irrotational* parts of flow velocity (3.3), as below ($\beta = const \neq 0$):

- for component $w_x$ $(x,y,z, t)$:

$$\frac{\partial w_x}{\partial t} = \nu \cdot \nabla^2 w_x \qquad (5.1)$$

- for component $w_y$ $(x,y,z, t)$:

$$\frac{\partial w_y}{\partial t} = \nu \cdot \nabla^2 w_y \qquad (5.2)$$

- for component $w_z$ $(x,y,z, t)$:

$$\frac{\partial w_z}{\partial t} = \nu \cdot \nabla^2 w_z \qquad (5.3)$$

Let us note that solutions of equations (5.1)-(5.3) with restrictions (4.3)-(4.4) should determine the appropriate *curvilinear* dynamical surface (time-dependent) which is bounding the subspace of existence of the proper exact solution of the Navier-Stokes equations (1.1)-(1.2). It means that one of the variables $\{x,y,z, t\}$ is to be depending on others, in this case.

So, the 3D non-stationary case of *partial* solution (3.3), (5.1)-(5.3) of Navier-Stokes equations (1.1)-(1.2) is to be reduced to the appropriate 2D-case (time-dependent).

## 6. <u>Discussion.</u>

To avoid ambiguity, we should clarify the appropriate procedure of finding the scale function $\gamma\,(x,\,y,\,z)$ in the expressions for presentation of a solution (3.3). Each of expression for partial derivative of function $\gamma\,(x,\,y,\,z)$ in regard to the chosen variable $\{x,\,y,\,z\}$ is included to the 3 of 4 equations as shown below (for example, $\partial\gamma/\partial x$):

$$U = -\gamma\cdot\left(\frac{2a}{1+(a^2+b^2)}\right),\quad V = -\gamma\cdot\left(\frac{2b}{1+(a^2+b^2)}\right),\quad W = \gamma\cdot\left(\frac{1-(a^2+b^2)}{1+(a^2+b^2)}\right),$$

(6.1)

$$1)\ \frac{\partial U}{\partial x}+\frac{\partial V}{\partial y}+\frac{\partial W}{\partial z}=0,\quad 2)\ \frac{\partial W}{\partial y}-\frac{\partial V}{\partial z}=0,\quad 3)\ \frac{\partial U}{\partial z}-\frac{\partial W}{\partial x}=0,\quad 4)\ \frac{\partial V}{\partial x}-\frac{\partial U}{\partial y}=0\ .$$

So, we should consider such 3 PDEs in (6.1), which are including $\partial\gamma/\partial x$ (for example), then we should eliminate the expression $\partial\gamma/\partial x$ by linear transformations from 2 of 3 of them; as a result, we should obtain 2 equations with linear expressions for $\partial\gamma/\partial y$ and $\partial\gamma/\partial z$.

Then we should repeat the procedure above and exclude also (by one linear transformation) the next expression for $\partial\gamma/\partial y$ from 1 of them:

$$\gamma \cdot \left( \frac{a}{b} \cdot \frac{\partial \left( \frac{1-(a^2+b^2)}{1+(a^2+b^2)} \right)}{\partial x} - \left( \frac{1-(a^2+b^2)}{2b} \right) \cdot \frac{\partial \left( \frac{2a}{1+(a^2+b^2)} \right)}{\partial x} \right) + \gamma \cdot \left( \frac{b}{a} \cdot \frac{\partial \left( \frac{1-(a^2+b^2)}{1+(a^2+b^2)} \right)}{\partial x} - \left( \frac{1-(a^2+b^2)}{2a} \right) \cdot \frac{\partial \left( \frac{2b}{1+(a^2+b^2)} \right)}{\partial x} \right) -$$

(6.2)

$$- \left( \frac{1-(a^2+b^2)}{2b} \right) \cdot \gamma \cdot \frac{\partial \left( \frac{2b}{1+(a^2+b^2)} \right)}{\partial y} + \left( \frac{1-(a^2+b^2)}{2a} \right) \cdot \gamma \cdot \frac{\partial \left( \frac{2a}{1+(a^2+b^2)} \right)}{\partial y} + \left( \frac{1-(a^2+b^2)}{2b} \right) \cdot \frac{\partial W}{\partial z} - \left( \frac{a^2+b^2}{ab} \right) \cdot \frac{\partial U}{\partial z} = 0 \; .$$

Thus, finally we should obtain from Eq. (6.2) the appropriate differential equation of the 1-st order for $\partial \gamma / \partial z$ in dependence on the function $\gamma(x,y,z)$, on the given functions $a$, $b$ (and their derivatives) – well, it is ready for integration as ODE in regard to $z$ variable! {*let us note that the final expression for $\gamma(x,y,z)$ should include also the arbitrary unknown functions, depending on variables $\{x,y\}$: $\gamma(x)_1 + \gamma(y)_2$*}:

$$H \cdot \frac{\partial \gamma}{\partial z} = G \cdot \gamma, \; H = \left( \frac{1-(a^2+b^2)}{2b} \right) \left( \frac{1-(a^2+b^2)}{1+(a^2+b^2)} \right) + \left( \frac{a^2+b^2}{ab} \right) \left( \frac{2a}{1+(a^2+b^2)} \right) = \left( \frac{1+a^2+b^2}{2b} \right)$$

$$G = \left( -\frac{a}{b} \cdot \frac{\partial \left( \frac{1-(a^2+b^2)}{1+(a^2+b^2)} \right)}{\partial x} + \left( \frac{1-(a^2+b^2)}{2b} \right) \cdot \frac{\partial \left( \frac{2a}{1+(a^2+b^2)} \right)}{\partial x} \right) +$$

(6.3)

$$+ \left( -\frac{b}{a} \cdot \frac{\partial \left( \frac{1-(a^2+b^2)}{1+(a^2+b^2)} \right)}{\partial x} + \left( \frac{1-(a^2+b^2)}{2a} \right) \cdot \frac{\partial \left( \frac{2b}{1+(a^2+b^2)} \right)}{\partial x} \right) +$$

$$+ \left( \frac{1-(a^2+b^2)}{2b} \right) \cdot \frac{\partial \left( \frac{2b}{1+(a^2+b^2)} \right)}{\partial y} - \left( \frac{1-(a^2+b^2)}{2a} \right) \cdot \frac{\partial \left( \frac{2a}{1+(a^2+b^2)} \right)}{\partial y} -$$

$$- \left( \frac{1-(a^2+b^2)}{2b} \right) \cdot \frac{\partial \left( \frac{1-(a^2+b^2)}{1+(a^2+b^2)} \right)}{\partial z} - \left( \frac{a^2+b^2}{ab} \right) \cdot \frac{\partial \left( \frac{2a}{1+(a^2+b^2)} \right)}{\partial z} \; .$$

We should repeat the procedure above (6.2)-(6.3) for each of the expressions $\partial \gamma/\partial y$ and $\partial \gamma/\partial z$ (for various chains of triplets of initial 4 equations). As a result, we should use all the 4 equations for such a procedure (where functions *a*, *b* are considered at any fixed meaning of *t*), but finally we should obtain all the parts for the expression of function $\gamma(x,y,z) = \gamma(x)_1 + \gamma(y)_2 + \gamma(z)_3$.

Let us note that at any meaning of time-parameter *t*, equations of a type (6.3) for the components $\gamma(x)_1$, $\gamma(y)_2$, $\gamma(z)_3$ of function $\gamma(x,y,z)$ should *spatially* bind (non-linearly) the components (5.1)-(5.3) of the *solenoidal* part of solution $w_x(x,y,z,t)$, $w_y(x,y,z,t)$ and $w_z(x,y,z,t)$. So, such a non-linear dependence should determine the appropriate *curvilinear* dynamical surface (time-dependent) which is bounding the subspace of existence of the proper exact solution of the Navier-Stokes equations (1.1)-(1.2). It means that one of the variables $\{x,y,z,t\}$ is to be depending on others, for each of components $w_x$, $w_y$ or $w_z$; thus, the *solenoidal* part of the solution is reduced to the 2D case (time-dependent).

Let us consider the example of exact solution (4.10)-(4.12) satisfying to Eqs. (6.1):

$$U = -\gamma \cdot \left(\frac{2a}{1+(\alpha^2+1)\cdot a^2}\right), \quad V = -\gamma \cdot \left(\frac{2\alpha \cdot a}{1+(\alpha^2+1)\cdot a^2}\right), \quad W = \gamma \cdot \left(\frac{1-(\alpha^2+1)\cdot a^2}{1+(\alpha^2+1)\cdot a^2}\right),$$

1) $-\dfrac{\partial \gamma}{\partial x}\left(\dfrac{2a}{1+(\alpha^2+1)\cdot a^2}\right) - \gamma \cdot \dfrac{\partial}{\partial x}\left(\dfrac{2a}{1+(\alpha^2+1)\cdot a^2}\right) - \alpha \cdot \dfrac{\partial \gamma}{\partial y}\left(\dfrac{2a}{1+(\alpha^2+1)\cdot a^2}\right) - \alpha \cdot \gamma \cdot \dfrac{\partial}{\partial y}\left(\dfrac{2a}{1+(\alpha^2+1)\cdot a^2}\right) +$

$+\dfrac{\partial \gamma}{\partial z}\left(\dfrac{1-(\alpha^2+1)\cdot a^2}{1+(\alpha^2+1)\cdot a^2}\right) + \gamma \cdot \dfrac{\partial}{\partial z}\left(\dfrac{1-(\alpha^2+1)\cdot a^2}{1+(\alpha^2+1)\cdot a^2}\right) = 0,$

(6.4)

2) $\dfrac{\partial \gamma}{\partial y}\left(\dfrac{1-(\alpha^2+1)\cdot a^2}{1+(\alpha^2+1)\cdot a^2}\right) + \gamma \cdot \dfrac{\partial}{\partial y}\left(\dfrac{1-(\alpha^2+1)\cdot a^2}{1+(\alpha^2+1)\cdot a^2}\right) + \alpha \cdot \dfrac{\partial \gamma}{\partial z}\left(\dfrac{2a}{1+(\alpha^2+1)\cdot a^2}\right) + \alpha \cdot \gamma \cdot \dfrac{\partial}{\partial z}\left(\dfrac{2a}{1+(\alpha^2+1)\cdot a^2}\right) = 0,$

3) $-\dfrac{\partial \gamma}{\partial z}\left(\dfrac{2a}{1+(\alpha^2+1)\cdot a^2}\right) - \gamma \cdot \dfrac{\partial}{\partial z}\left(\dfrac{2a}{1+(\alpha^2+1)\cdot a^2}\right) - \dfrac{\partial \gamma}{\partial x}\left(\dfrac{1-(\alpha^2+1)\cdot a^2}{1+(\alpha^2+1)\cdot a^2}\right) - \gamma \cdot \dfrac{\partial}{\partial x}\left(\dfrac{1-(\alpha^2+1)\cdot a^2}{1+(\alpha^2+1)\cdot a^2}\right) = 0,$

4) $-\alpha \cdot \dfrac{\partial \gamma}{\partial x}\left(\dfrac{2a}{1+(\alpha^2+1)\cdot a^2}\right) - \alpha \cdot \gamma \cdot \dfrac{\partial}{\partial x}\left(\dfrac{2a}{1+(\alpha^2+1)\cdot a^2}\right) + \dfrac{\partial \gamma}{\partial y}\left(\dfrac{2a}{1+(\alpha^2+1)\cdot a^2}\right) + \gamma \cdot \dfrac{\partial}{\partial y}\left(\dfrac{2a}{1+(\alpha^2+1)\cdot a^2}\right) = 0.$

At 1-st step, we should multiply Eq. 1) of Eqs. (6.4) on coefficient (- α), then we should sum it to the Eq. 4); as a result, we obtain:

$$(\alpha^2+1)\cdot\frac{\partial \gamma}{\partial y}\cdot 2a + (\alpha^2+1)\cdot\gamma\cdot\frac{2(1-(\alpha^2+1)\cdot a^2)}{(1+(\alpha^2+1)\cdot a^2)}\cdot\frac{\partial a}{\partial y} - \alpha\frac{\partial \gamma}{\partial z}\cdot(1-(\alpha^2+1)\cdot a^2) + \frac{4(\alpha^2+1)\cdot a\cdot\alpha\,\gamma}{(1+(\alpha^2+1)\cdot a^2)}\cdot\frac{\partial a}{\partial z} = 0,$$

- but Eq. 2) could be transformed as below:

$$\frac{\partial \gamma}{\partial y}\cdot(1-(\alpha^2+1)\cdot a^2) - \gamma\cdot\frac{4(\alpha^2+1)\cdot a}{(1+(\alpha^2+1)\cdot a^2)}\cdot\frac{\partial a}{\partial y} + 2\alpha\cdot\frac{\partial \gamma}{\partial z}\cdot a + 2\alpha\cdot\gamma\cdot\frac{(1-(\alpha^2+1)\cdot a^2)}{(1+(\alpha^2+1)\cdot a^2)}\cdot\frac{\partial a}{\partial z} = 0 \;.$$

Last two equations let us determine the appropriate expressions for $\partial a/\partial y$ and $\partial a/\partial z$:

$$\frac{\partial a}{\partial y} = \left(\frac{\alpha}{2(\alpha^2+1)\cdot\gamma}\right)\cdot(1+(\alpha^2+1)\cdot a^2)\cdot\frac{\partial \gamma}{\partial z}, \quad \frac{\partial a}{\partial z} = -\left(\frac{1}{2\alpha\cdot\gamma}\right)\cdot(1+(\alpha^2+1)\cdot a^2)\cdot\frac{\partial \gamma}{\partial y}, \qquad (6.5)$$

Besides, Eq. 3) let us obtain the appropriate expression for $\partial a/\partial x$ as below:

$$\frac{\partial a}{\partial x} = \left((1-(\alpha^2+1)\cdot a^2)\cdot\frac{\partial \gamma}{\partial x} - \frac{(1-(\alpha^2+1)\cdot a^2)}{\alpha}\cdot\frac{\partial \gamma}{\partial y} + 2a\cdot\frac{\partial \gamma}{\partial z}\right)\cdot\frac{(1+(\alpha^2+1)\cdot a^2)}{4(\alpha^2+1)\cdot a\cdot\gamma} \qquad (6.6)$$

Let us consider a proper simplifying assumptions of the expression (6.6) for $\partial a/\partial x$ (here below we assume $\{\partial \gamma/\partial x, \partial \gamma/\partial y\} \neq 0$)

$$\frac{\partial \gamma}{\partial z} = 0, \quad \frac{\partial \gamma}{\partial x} = \frac{1}{\alpha}\cdot\frac{\partial \gamma}{\partial y}, \qquad (6.7)$$

- which should transform (6.5)-(6.6) as below

$$\frac{\partial a}{\partial x} = 0, \quad \frac{\partial a}{\partial y} = 0, \quad \Rightarrow \quad w_y = w_y(z,t) \qquad (6.8)$$

We should additionally note that one of the partial (possible) solutions for Eqs. (6.7) is proved to be $\gamma = \gamma_0\cdot\exp(-kx - \alpha\cdot ky)$, where $k, \gamma_0 = \text{const}$.

So, taking into consideration the expression for $a(x,y,z, t)$ from (4.12) and a proper expression for $\partial a/\partial z$ from Eq. (6.5), we obtain as below

$$a(x, y, z, t) = \frac{\tan\left(\frac{(\alpha^2+1)}{2} \cdot \int w_y \, dt\right)}{(\alpha^2+1)}, \quad \frac{\partial a}{\partial z} = -\left(\frac{1}{2\alpha \cdot \gamma}\right) \cdot (1 + (\alpha^2+1) \cdot a^2) \cdot \frac{\partial \gamma}{\partial y},$$

$$\Rightarrow \int\left(\frac{\partial w_y}{\partial z}\right) dt = -\left(\frac{1}{\alpha \cdot (\alpha^2+1) \cdot \gamma}\right) \cdot \frac{\partial \gamma}{\partial y} \cdot \left(\alpha^2 \cdot \cos^2\left(\frac{(\alpha^2+1)}{2} \cdot \int w_y \, dt\right) + 1\right) \quad (6.9)$$

If we differentiate both the parts of last equation (6.9) in regard to the time-parameter $t$, it should yield

$$\frac{\partial w_y}{\partial z} = \left(\frac{\alpha}{2\gamma}\right) \cdot \frac{\partial \gamma}{\partial y} \cdot \sin\left((\alpha^2+1) \cdot \int w_y \, dt\right) \cdot w_y, \quad (6.10)$$

We could transform (6.10) as below:

$$(\alpha^2+1) \cdot \int w_y \, dt = \arcsin\left(\frac{2\gamma \cdot \frac{\partial w_y}{\partial z}}{\alpha \cdot \frac{\partial \gamma}{\partial y} \cdot w_y}\right) \quad (6.11)$$

Thus, finally we obtain from (6.11) and (4.12)

$$a(x, y, z, t) = \frac{\tan\left(\frac{1}{2}\arcsin\left(\frac{2\gamma \cdot \frac{\partial w_y}{\partial z}}{\alpha \cdot \frac{\partial \gamma}{\partial y} \cdot w_y}\right)\right)}{(\alpha^2+1)}, \quad (6.12)$$

- where $w_y(z, t)$ – is the *non-stationary* solution of one-dimensional Eq. (5.2):

$$\frac{\partial w_y}{\partial t} = \nu \cdot \nabla^2 w_y$$

## 7. Final presentation of exact solution.

Let us present exact non-stationary solution $\{p, \boldsymbol{u}\}$ of Navier-Stokes equations (1.1)-(1.2) in a final form ($\beta = const \neq 0$):

$$\frac{\nabla p}{\rho} = -\nabla \phi - \frac{1}{2}\nabla\{(\vec{u}_p + \vec{u}_w)^2\}, \quad \vec{u} \equiv \vec{u}_p + \vec{u}_w, \quad \{\nabla \cdot \vec{u}_w \equiv 0, \quad \nabla \times (\vec{u}_p) \equiv 0\},$$

$$\vec{u}_p \equiv \{U, V, W\}, \quad \vec{u}_w \equiv \{\beta \cdot w_x, \beta \cdot w_y, \beta \cdot w_z\}, \tag{7.1}$$

$$U = -\gamma \cdot \left(\frac{2a}{1+(a^2+b^2)}\right), \quad V = -\gamma \cdot \left(\frac{2b}{1+(a^2+b^2)}\right), \quad W = \gamma \cdot \left(\frac{1-(a^2+b^2)}{1+(a^2+b^2)}\right),$$

- where $\rho$ is the fluid density, $\phi$ – is the potential of external force, acting on a fluid.

Besides, the components $a$, $b$ of the *irrotational* parts of flow velocity (3.3) are the solutions of system (3.5), see examples: - 1) Eq. (4.3); - 2) Eq. (4.4); - 3) Eqs. (4.9), (4.6); and - 4) Eqs. (4.12), (4.10). Let us note that function $\gamma(x,y,z)$ along with the functions $\{U, V, W\}$ should satisfy to the equations (2.6)-(2.7) for such an expressions.

Note also that $w_i$ in (7.1), $i = \{x, y, z\}$, are the appropriate solutions of the equations (5.1)-(5.3) (under the proper initial conditions). If we assume *self-similarity* ansatz [6] $w_i = w_i(t) \cdot w_i(x,y,z)$, it should yield from (5.1)-(5.3) as below

$$w_i = \exp(-\alpha_i^2 \cdot t) \cdot w_i(x, y, z), \tag{7.2}$$

$$\nabla^2 w_i(x, y, z) + \left(\frac{\alpha_i^2}{\nu}\right) \cdot w_i(x, y, z) = 0, \tag{7.3}$$

- where $\alpha_i = const$; the last equation (7.3) is known to be the *spatial* Helmholtz equation [9], which describes evolution of spatial part $w_i(x,y,z)$ of the component $w_i$.

As for the examples of solution of *spatial* Helmholtz (7.3) equation, let us remember well-known spherical-wave solutions and their non-symmetric generalization [9].

## 8. Conclusion.

In fluid mechanics, a lot of authors have been executing their researches to obtain the analytical solutions of Navier-Stokes equations [10], even for 3D case of *compressible* gas flow [11] or for 3D case of *non-stationary* flow of incompressible fluid [12]. But there is an essential deficiency of non-stationary solutions indeed.

In our presentation, we explore the ansatz of derivation of *non-stationary* solution for the Navier-Stokes equations in the case of incompressible flow, which was suggested previously.

In general case, such a solution should be obtained from the mixed system of 2 *Riccati* ordinary differential equations (in regard to time-parameter $t$). But we find an elegant way to simplify it to the proper analytical presentation of exact solution (such a solution is exponentially decreasing to zero for $t$ going to infinity $\infty$). Also it has to be specified that the solutions that are constructed can be considered as a class of perturbation absorbed exponentially as $t$ going to infinity $\infty$ by the null solution [13].

The uniqueness of the presented solutions is not considered. In this respect we confine ourselves to mention the paper [14], in which all the difficulties concerning the uniqueness in unbounded domain are remarked.

## Acknowledgements

I am thankful to unknown esteemed Reviewer for valuable comprehensive advices in preparing of this manuscript.

I devote this article to my darling wife who is preparing for giving birth to our beloved daughter.

[13]. Berker, R. (1963), *Intégration des équations du mouvement d'un fluide visqueux incompressible.* In: Encyclopedia of Physics (Flügge, S., ed.) Vol. 8, no. 2. Berlin-Göttingen-Heidelberg: Springer 1963.

[14]. Galdi G.P., Rionero S. (1979), *The weight function approach to uniqueness of viscous flows in unbounded domains*, Archive for Rational Mechanics and Analysis, Volume 69, Issue 1, pp 37-52.